\documentclass[12pt,a4paper]{article}
 \usepackage[top=40mm,bottom=40mm]{geometry}

\usepackage{amsmath,defns,reducecode,eulervm,natbib,url,graphicx,microtype}

\usepackage{hyperref} 
\hypersetup{pdftex,
            backref=true,
            hyperindex=true,
            colorlinks=true,
            citecolor=blue,
            bookmarks=true,
            breaklinks=true}
            
\allowdisplaybreaks

\newcommand{\uu}{{\bar u}}
\newcommand{\vv}{{\bar v}}
\newcommand{\cc}{{\bar c}}
\newcommand{\bq}{{\bar q}}

\newcommand{\ros}{{\dot\varepsilon}}
\renewcommand{\vec}[1]{\text{\boldmath$#1$}}

\allowdisplaybreaks

\title{Modelling suspended sediment in environmental turbulent fluids}

\author{Meng Cao
\thanks{School of Mathematical Sciences,
University of Adelaide, South Australia 5005.  \url{mailto:meng.cao@adelaide.edu.au} or  \url{mailto:mengcao1188216@gmail.com}}
\and 
A.~J. Roberts
\thanks{School of Mathematical Sciences,
University of Adelaide, South Australia 5005.  \url{mailto:anthony.roberts@adelaide.edu.au}}
}

\begin{document}
    
\maketitle

\begin{abstract}
Modelling sediment transport in environmental turbulent fluids is a challenge.
This article develops a sound model of the lateral transport of suspended sediment in environmental fluid flows such as floods and tsunamis.
The model is systematically derived from a 3D turbulence model based on the Smagorinski large eddy closure.
Embedding the physical dynamics into a family of problems and analysing linear dynamics of the system, centre manifold theory indicates the existence of slow manifold parametrised by macroscale variables.
Computer algebra then constructs the slow manifold in terms of fluid depth, depth-averaged lateral velocities, and suspended sediment concentration.
The model includes the effects of sediment erosion, advection, dispersion, and also the interactions between the sediment and turbulent fluid flow.
Vertical distributions of the velocity and concentration in steady flow agree with the established experimental data.
Numerical simulations of the suspended sediment under large waves  show that the developed model predicts physically reasonable phenomena.
\end{abstract}

\tableofcontents

\section{Introduction}
\label{sec:intro}

\begin{figure}
\centering
\begin{tabular}{c@{}c}
\includegraphics[width=0.8\textwidth]{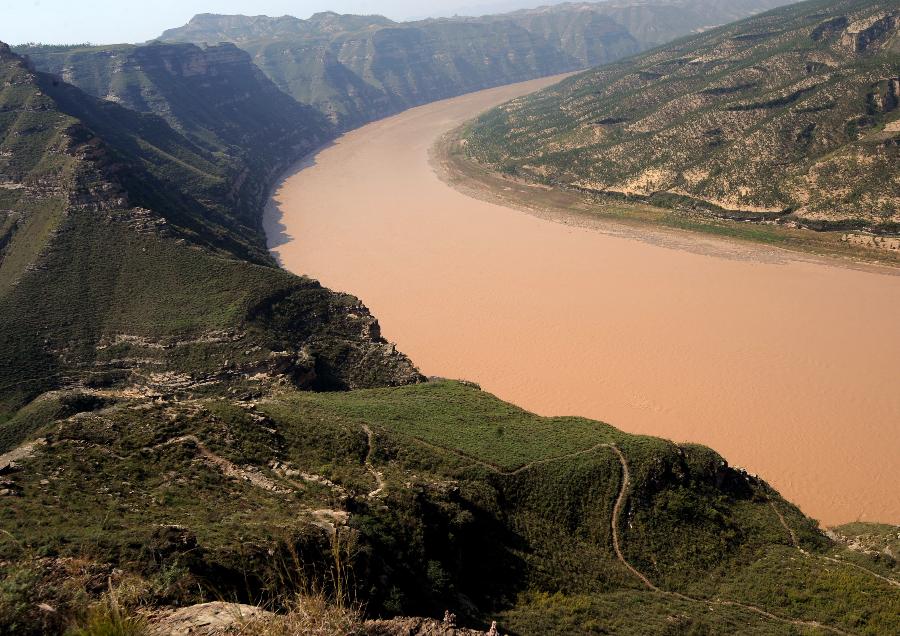}
\end{tabular}
\caption{Scene of Yellow River turning in Shilou, China's Shanxi~(\url{http://news.xinhuanet.com}).
This river carries vast amounts of suspended sediment.
}
\label{IntroYR}
\end{figure}

Environmental turbulent fluids, such as rivers, floods and tsunamis, always carry amounts of sediment. 
For example, Figure~\ref{IntroYR} shows the Yellow River in China which is famous for carrying large amounts of sediment.
Modelling the sediment in these environmental fluids is important for studying and predicting changes of the morphology and topography.
We aim to develop a model to appropriately model the suspended sediment in turbulent fluid flows via systematic resolution of the physical processes.

Our modelling, which is based on dynamical systems theory instead of conventional depth-averaging, resolves out-of-equilibrium interactions between the varying turbulence and suspended sediment.
Most previous work studied suspended sediment in uniform flows~\cite[e.g.]{Hunt:1954fk, Rijn:1984fk, Celik:1988fk, Deigaard1992} or by depth-averaging flow and sediment equations~\cite[e.g.]{Wu:2000uq, Pittaluga:2003fk}. 
We explore the implications of changing the theoretical base from depth-averaging to a slow manifold of the turbulent Smagorinski large eddy closure.

\cite{Cao:2012fk} initially developed a 2D lateral model for environmental fluid flows, derived from a 3D turbulence model based on Smagorinski large eddy closure.
This model includes the effects and interactions of inertia, advection, bed drag, gravitational forcing and turbulent dissipation with minimal assumptions.
The innovation here is that the turbulent modelling and dynamics includes and interacts with the suspended sediment transport \cite[]{Cao2014}.
A slow manifold is found for the out-of-equilibrium dynamics of the coupled turbulent sediment system.
We choose to parametrise the slow manifold in terms of emergent depth-averaged quantities.
The evolution of these depth-averaged quantities on the slow manifold governs the dynamics of the suspended sediment in the turbulent fluid flows.

Consider a turbulent flow of depth $h(x,y,t)$ flowing along a bed $z=b(x,y)$ with a mean slope~$\tan\theta$ and carrying sediment.
Sections~\ref{sec:detailEq}--\ref{sec:detailBc} detail equations of the Reynolds-averaged Navier--Stokes \pde{}s, the advection-diffusion equation and boundary conditions on the free surface and the mean bed. 
Section~\ref{smag:Tbc} uses the dynamical systems theory of centre manifolds~\cite[e.g.]{Roberts1988, Potzsche:2006uq} to analyse the governing equations and derive the following non-dimensional suspended sediment model of the horizontal evolution of the depth~$h(x,y,t)$, the lateral depth-averaged velocities~$\uu(x,y,t)$ and~$\vv(x,y,t)$ and depth-averaged concentration~$\cc(x,y,t)$:
\begin{subequations}\label{eqs:Intro}
\begin{align}
\D{t}{h}\approx&-\D{x}{}\left(h\uu\right)-\D{y}{}\left(h\vv\right)\,,\label{sed:IntroH}
\\
\D{t}{\bar u}\approx&{}-0.00293\frac{\bar u\bq}{h}+0.993\left[\tan\theta-\D{x}{}\left(h+b\right)\right]
-1.025\uu\D x\uu-1.017\vv\D y\uu
\nonumber\\&{}
-0.298(s-1)h\D x\cc
\,,\label{sed:IntroU}
\\
\D{t}{\vv}\approx&{}-0.00293 \frac{\vv\bq}{h}-0.993\D{y}{}\left(h+b\right)
-1.025\vv\D y\vv-1.017\uu\D x\vv
\nonumber\\&{}
-0.298(s-1)h\D y\cc
\,,\label{sed:IntroV}
\\
\D t\cc\approx{}&{}
-\frac{w_f}{h}\left(0.938\cc -0.984c_{ae}\right)
-1.007\exp\left(-3.073\frac{w_f}{\bq}\right)\left(\uu\D x\cc+\vv\D y\cc\right),\label{sed:IntroC}
\end{align}
\label{Intro:model}%
\end{subequations}
where $\bq=\sqrt{\uu^2+\vv^2}$  the mean local flow speed, constant~$w_f$ is the falling velocity of the sediment, and constant~$c_{ae}$ is an equilibrium reference concentration on the mean bed~$z=b$.
The effective momentum equations~\eqref{sed:IntroU}--\eqref{sed:IntroV}, and the more refined version~\eqref{smag:u}--\eqref{smag:v}, include the effects of gravitational forcing, bed drag, self-advection,  turbulent dissipation, and sediment induced flow.
The sediment \pde~\eqref{sed:IntroC} includes the sediment erosion and deposition, and the advection.
Although this model is expressed in terms of depth-averaged lateral velocities and concentration, they are derived not by depth-averaging, but instead by systematically accounting for interaction between vertical profiles and horizontal gradients of the velocity, the concentration, the stress, bed drag, lateral space variations and bed topography.

That the systematically reduced model~\eqref{eqs:Intro} has many terms so close to established models is a testament to the robustness of conservation of mass and momentum principles that underly traditional derivations.
Additionally, our dynamical systems approach resolves finer microscale and out-of-equilibrium interactions that are active in more physically complicated systems. 

Linear analysis in section~\ref{smag:Tlinear} finds the spectrum  supports the existence of a slow manifold.
The computer algebra of Appendix~\ref{caSmag} then constructs the slow manifold of the system.
Then section~\ref{sec:model} derives the reduced model~\eqref{eqs:Intro} of the turbulent flow and suspended sediment on the slow manifold.

Section~\ref{sec:Ver} discusses the predicted vertical distribution of the velocity and concentration fields in steady flow.
Agreement with established experimental data~\cite[e.g.]{Schultz:2007uq, Schultz:2013fk, Celik:1988fk} indicates that our model is reasonable.
Section~\ref{sec:Num} numerically explores the suspended sediment under large waves by the comprehensive suspended sediment model~\eqref{sed:cH} coupled with the turbulence model~\eqref{smag:h}--\eqref{smag:v}.
The numerical results indicate that our model reasonably describes the dynamics of the suspended sediment in turbulent flows.

\section{Construct the sediment transport model}
\label{sec:detail}

This section describes the derivation of the reduced models for turbulent flow and suspended sediment.
First, section~\ref{sec:detailEq} details the 3D continuity, Navier--Stokes and advection-diffusion equations of the turbulent fluid flows and suspended sediment, whereas section~\ref{sec:detailBc} records the boundary conditions of the flow and the sediment on the free surface and the mean bed.
Second, section~\ref{smag:Tbc} embeds these equations in a family of equations with modified tangential stresses on the free surface to establish the existence of an appropriate slow manifold. 
The linear analysis of section~\ref{smag:Tlinear} supports the emergence of the slow manifold from the dynamics in the system.
The computer algebra of Appendix~\ref{caSmag} handles the details of the  construction of the slow manifold model that is summarised in section~\ref{sec:model}.

\subsection{The governing equations of the flow and sediment}
\label{sec:detailEq}

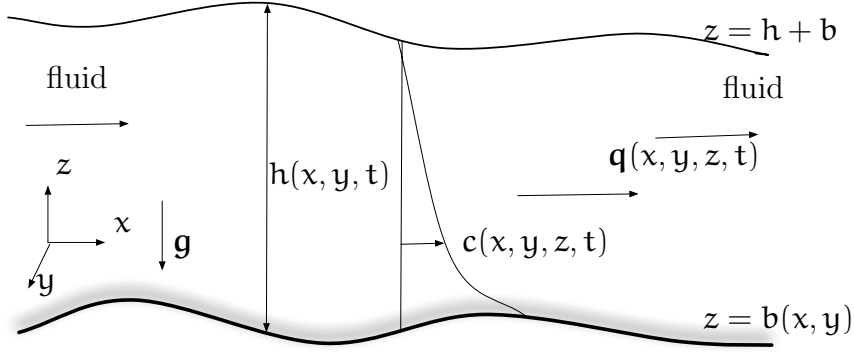
\begin{figure}
\begin{center}
\setlength{\unitlength}{0.7ex}
\begin{picture}(120,50)(0,0)
\put(10,0){\includegraphics[width=86\unitlength]{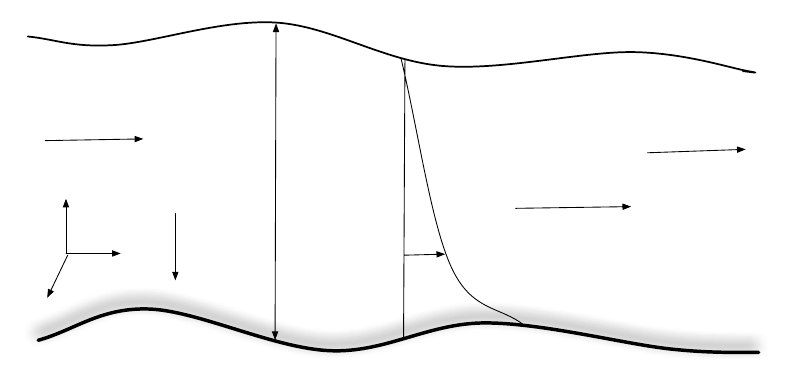}}
\put(24,15){$x$} \put(16,9){$y$} \put(18,21){$z$}
\put(17,30){{fluid}} \put(87,29){{fluid}}
\put(40,20){$h(x,y,t)$} \put(75,22){$\vec q(x,y,z,t)$}
\put(30,13){$\vec g$}
\put(60,13){$c(x,y,z,t)$}
\put(85,5){$z=b(x,y)$} \put(85,35){$z=h+b$}
\end{picture}
\end{center}
\caption{This diagram depicts the suspended sediment in turbulent flow with~$(x,y,z)$ coordinate system.
The fluid of depth~$h(x,y,t)$ flows down the sloped bed at the turbulent mean velocity~$\vec q(x,y,z,t)$.
The turbulent mean concentration of the suspended sediment is~$c(x,y,z,t)$.
Denote the mean bed $z=b(x,y)$, the free surface~$z=h+b$\,, and the gravity is~$g$.
}
\label{sedDiagram}
\end{figure}

Consider the turbulent flow flowing along a bed carrying sediment.
This work only considers the suspended sediment and neglects the bed load transport on the mean bed.
Figure~\ref{sedDiagram} depicts the diagram of the suspended sediment in the turbulent flow.
Define a coordinate system with~$x,y$ for the lateral directions,  and~$z$ for the direction normal to the mean slope.

The fluid of depth~$h(x,y,t)$ flows down the sloped mean bed~$z=b(x,y)$ at the turbulent mean velocity~$\vec q(x,y,z,t)$; the velocity vector $\vec q=(u,v,w)$ in the $(x,y,z)$~directions, respectively.
The term `mean bed' refers to the smooth average bed over an ensemble of turbulence and bed roughness realisations: the physical bed in any one realisation is envisaged to be rough, just like the physical fluid velocity field would have rapid spatial variations in any one realisation.
Denote the turbulent mean pressure field by~$p(x,y,z,t)$.
The suspended sediment has a turbulent mean concentration~$c(x,y,z,t)$~(volume fraction).
The mean bed $z=b(x,y)$ has an overall slope~$\tan\theta$ in the $x$-direction.

We assume the particles of the suspended sediment have small sizes, and all the particles of the suspended sediment have the same falling velocity, namely that of a sphere of diameter~$d$.

The nondimensional governing partial differential equations for the incompressible, three dimensional, turbulent mean fluid fields are the Reynolds-averaged continuity and momentum equations, 
\begin{align}&
\divv\vec q=\D xu+\D yv+\D zw=0\,,\label{sec:mass}
\\&
\D{t}{\vec q}+\vec q\cdot\grad\vec q=-\grad \vec p+\divv\tau+\vec g\,,
\label{sec:mom}
\end{align}
and for the suspended sediment is the advection-diffusion equation,
\begin{equation}
\D tc+\divv(\vec qc)=-\divv(w_fc\vec n_g)+\divv(\nu\grad c)\,,
\label{sec:con}
\end{equation}
where $w_f$~is the falling velocity of the sediment particles, $\nu$~is the eddy viscosity, and the vector $\vec g=(\tan\theta,0,-1)$ is the direction of the nondimensional gravity.

The suspended sediment influences the fluid turbulence.
We assume the mixing density~$\rho_{\text{mix}}$ of the fluid and suspended sediment satisfies
\begin{align}
\frac{1}{\rho_{\text{mix}}}&=\frac{1}{\rho+c(\rho_m-\rho)}=\frac{1}{\rho}\frac{1}{1+c(s-1)}
\nonumber\\&
=\frac{1}{\rho}\left[1-c(s-1)+c^2(s-1)^2\right]+\mathcal O([c(s-1)]^3)\,,
\label{sed:rho}
\end{align}
where $\rho_m$~is the sediment density, $\rho$~the density of the fluid and $s=\rho_m/\rho$ is the relative density.

The \textsc{pde}s~\eqref{sec:mass}--\eqref{sec:con} are nondimensional and derived upon using the characteristic depth~$H$ of the turbulent fluid as the length scale, the long wave speed~$\sqrt{g_zH}$ as the velocity scale, and the fluid density~$\rho$ as the reference density.
Thus, the mixing density~$\rho_{\text{mix}}\approx1-c(s-1)$.

The variable~$\tau$ is the turbulent mean deviatoric stress tensor, which is approximated using the eddy viscosity~$\nu$. 
We use the Smagorinski eddy closure to approximate the turbulence stresses.
\cite{Cao:2012fk} and \cite{Georgiev2008} detailed this eddy closure and expressed the mean deviatoric stress tensor
\begin{equation}
\tau_{ij}=\nu\left( \D{x_j}{u_i} +\D{x_i}{u_j}\right)=c_th^2\ros\left( \D{x_j}{u_i} +\D{x_i}{u_j}\right)\,,
\label{smag:tau}
\end{equation}
where the magnitude of the second invariant of the strain-rate tensor satisfies $|\ros|^2=\sum_{i,j}\left( \D{x_j}{u_i} +\D{x_i}{u_j}\right)^2$.

The falling velocity~$w_f$ is related to the mean particle effective size~$d$, the relative density~$s$ and the gravity~$g$.
We set the falling velocity~\cite[e.g.]{Deigaard1992}
\begin{equation}
w_f=\sqrt{\frac{4(s-1)gd}{3c_D}}\,.
\label{sed:wf}
\end{equation}
The drag coefficient $c_D\approx1.4$ for the large grain Reynolds number of natural sands, typically~$\re>500$~\cite[e.g.]{Deigaard1992}.
Typically, a quartz sediment has a relative density $s=2.65$~\cite[\S7.1]{Chanson:2004fk}.

\subsection{The boundary conditions of the flow and sediment}
\label{sec:detailBc}

We formulate boundary conditions on the mean bed $z=b(x,y)$ and on the free surface $z=\eta(x,y,t)=h(x,y,t)+b(x,y)$ in terms of the turbulent mean velocity field~$\vec q(x,y,z,t)$, turbulent mean concentration~$c(x,y,z,t)$, the fluid depth~$h(x,y,t)$, and the turbulent mean pressure~$p(x,y,z,t)$. 

First formulate the boundary conditions for the turbulent flows.
\begin{itemize}
\item On the mean bed, no fluid penetrating the ground requires
\footnote{This and the following boundary conditions are expressed in terms of ensemble mean quantities.  Consequently,  terms in the mean of the products of fluctuations might appear and a closure for them invoked \cite[\S2.2.6]{Cao2014}.  We assume the closure that such products of fluctuations are negligible in the boundary conditions.}
  \begin{equation}
  w=ub_x+vb_y \quad\text{on } z=b\,.
  \label{smag:nopenetrate}
  \end{equation}
\item Positing a slip law on the mean bed to account for a negligibly thin turbulent boundary layer gives
 \begin{align}&
 \frac{1}{\sqrt{1+b_x^2}}(u+wb_x)=\frac{c_uh}{\sqrt{1+b_x^2+b_y^2}}(u+wb_x)_{\vec n}
  \quad\text{on } z=b\,,\label{smag:slipu}\\&
 \frac{1}{\sqrt{1+b_y^2}}(v+wb_y)=\frac{c_uh}{\sqrt{1+b_x^2+b_y^2}}(v+wb_y)_{\vec n}
 \quad\text{on } z=b\,,\label{smag:slipv}
 \end{align}
where the derivative~$\partial_{\vec n}=-b_x\partial_x-b_y\partial_y+\partial_z$ and the constant $c_u\approx1.85$ matches open channel flow observations~\cite[e.g.]{Roberts2008}. 
\item On the free surface (that is, on its turbulent mean position), the kinematic condition that no fluid crosses the free surface is 
 \begin{equation}
 \eta_t+u\eta_x+v\eta_y=w \quad\text{on } z=\eta=h+b\,.
 \label{smag:kinetic}
 \end{equation}
\item Zero turbulent mean stress normal to the free surface indicates that on \(z=\eta\)
 \begin{equation}
    -p+\frac{\tau_{33} -2\eta_x\tau_{13} -2\eta_y\tau_{23}
    +\eta_x^2\tau_{11} +2\eta_x\eta_y\tau_{12}+\eta_y^2\tau_{22}}
    {1+\eta_x^2+\eta_y^2}
     =0\,.
    \label{smag:ttz}
 \end{equation}
\item  No turbulent mean, tangential stresses at the free surface indicates that on \(z=\eta\)
\begin{eqnarray}&&
    (1-\eta_x^2)\tau_{13}+\eta_x(\tau_{33}-\tau_{11})-\eta_y(\tau_{12}+\eta_x\tau_{23})=0
\,,
    \label{smag:ttxN} \\&&
    (1-\eta_y^2)\tau_{23}+\eta_y(\tau_{33}-\tau_{22})
    -\eta_x(\tau_{12}+\eta_y\tau_{13})=0
\,.
    \label{smag:ttyN}
\end{eqnarray}
\end{itemize}
There are two boundary conditions for the suspended sediment.
\begin{itemize}
\item On the free surface, the sediment flux normal to the surface is zero, which requires 
\begin{equation}
(\eta_x\tan\theta+1)w_fc+\nu\left(-\eta_xc_x-\eta_yc_y+c_z\right)=0 \quad\text{on } z=\eta\,.
\label{sed:bcsurface}
\end{equation}
\item On the mean bed, the upward net flux across the mean bed comes from the entrainment due to the fluid turbulence, that is
 \begin{equation}
 \nu\left(b_xc_x+b_yc_y-c_z\right)=\left(b_x\tan\theta+1\right)w_fc_{ae} \quad\text{on } z=b\,,
 \label{sed:bcbedf}
 \end{equation}
where, following the work of~\cite{Rijn:1984fk}, the equilibrium reference concentration is approximated in terms of shear velocity~$u_f$ and mean particle size~$d$ as
 \begin{equation}
 c_{ae}\approx0.075\frac{u_f^3}{d^{0.8}}\,.
 \label{sed:cae1}
 \end{equation}
\end{itemize}

The nondimensional \textsc{pde}s~\eqref{sec:mass}--\eqref{sec:con}, together with boundary conditions~\eqref{smag:nopenetrate}--\eqref{sed:bcbedf} govern the dynamics of the turbulent flow and suspended sediment.



\subsection{Embed the physical problem in a family of problems}
\label{smag:Tbc}

In order to provide theoretical support for the model redction, we embed the surface conditions~\eqref{smag:ttxN}--\eqref{smag:ttyN} in a family of conditions that modify the tangential stresses to have an artificial forcing proportional to the square of the local, free surface, velocity:
\begin{eqnarray}&&
    (1-\eta_x^2)\tau_{13}+\eta_x(\tau_{33}-\tau_{11})
    -\eta_y(\tau_{12}+\eta_x\tau_{23})
    \nonumber\\&&{}
    = \frac{(1-\gamma)\sqrt2c_t}{(1+c_u)(1+2c_u)} u\sqrt{u^2+v^2}
    \quad\text{on } z=\eta=h+b\,,
    \label{smag:ttx}
    \\&&
    (1-\eta_y^2)\tau_{23}+\eta_y(\tau_{33}-\tau_{22})
    -\eta_x(\tau_{12}+\eta_y\tau_{13})
    \nonumber\\&&{}
    = \frac{(1-\gamma)\sqrt2c_t}{(1+c_u)(1+2c_u)} v\sqrt{u^2+v^2}
    \quad\text{on } z=\eta=h+b\,.
    \label{smag:tty}
\end{eqnarray}
When evaluated at parameter $\gamma=1$ these artificial right-hand side becomes zero and the artificial surface conditions~\eqref{smag:ttx}--\eqref{smag:tty} reduce to the physical surface condition~\eqref{smag:ttxN}--\eqref{smag:ttyN}.  

However, when the artificial parameter $\gamma=0$, and when the mean slope and the lateral derivatives are negligible ($\tan\theta=\partial_x=\partial_y=0$), then the boundary conditions~\eqref{smag:ttx}--\eqref{smag:tty} reduce to $\nu\D zu=\nu (u/\eta)$ and $\nu\D zv=\nu (v/\eta)$.
In this case, two neutral modes of the dynamics are the lateral shear flow~$(u,v)\propto z/h+c_u$\,.
\footnote{The Euler parameter of a toy problem suggests introducing a factor $(1-\gamma/6)$ into the left-hand side of the tangential stress boundary conditions~\eqref{smag:ttxN}--\eqref{smag:ttyN} in order to improve convergence in the parameter~$\gamma$ when evaluated at the physically relevant $\gamma=1$\,. 
For the moment, this work omits such a factor.}

Analogously we embed the sediment boundary conditions~\eqref{sed:bcsurface} and~\eqref{sed:bcbedf} in the family
\begin{align}
&
\nu\left(-\eta_xc_x-\eta_yc_y+\gamma_cc_z\right)+2(1-\gamma_c)\nu\frac{c}{h}
\nonumber\\&\quad
{}+\left[1+(1-\gamma_c)\frac{w_f}{6}\right]\left(1+\eta_x\tan\theta\right)w_fc=0
\quad\text{on}\quad
z=h+b\,,\label{sed:bcsurfaceG}
\\&
\nu\left[-b_xc_x-b_yc_y+(2-\gamma_c)c_z\right]+2(1-\gamma_c)\nu\frac{c}{h}
\nonumber\\&\quad
+\left[1+(1-\gamma_c)\frac{w_f}{6}\right]\left(1+b_x\tan\theta\right)w_fc_{ae}=0
\quad\text{on}\quad
z=b\,.\label{sed:bcbedfG}
\end{align}
Upon setting the embedded parameter~$\gamma_c=1$\,, the boundary conditions~\eqref{sed:bcsurfaceG}--\eqref{sed:bcbedfG} recover the original physical sediment boundary conditions~\eqref{sed:bcsurface}--\eqref{sed:bcbedf}.
When the embedded parameter $\gamma_c=0$\,, boundary conditions~\eqref{sed:bcsurfaceG}--\eqref{sed:bcbedfG} are part of an artificial problem which is used to find a slow manifold in the physical system.
The extra term~$w_f/6$ in equations \eqref{sed:bcsurfaceG}--\eqref{sed:bcbedfG} ensures conservation to errors~$\mathcal O(w_f^4)$.
Without such term, the model only conserves sediment to errors~$\mathcal O(w_f^2)$.
The reason for implementing such high order errors is that the information about different falling velocities comes into the model in~$\mathcal O(w_f^2)$ terms.

When the artificial parameter~$\gamma_c=0$, the lateral gravity and lateral derivatives are negligible ($\tan\theta=\partial_x=\partial_y=0$), and the falling velocity~$w_f=0$, boundary condition~\eqref{sed:bcsurfaceG} requires the concentration~$c=0$ on the free surface, and boundary condition~\eqref{sed:bcbedfG} indicates
\begin{equation*}
\D zc+\frac{c}{h}=0
\quad\text{on}\quad 
z=b\,.
\end{equation*}
Together with the sediment \pde~\eqref{sec:con}, these imply a neutral mode of the sediment dynamics is $c\propto 1-z/h$ when artificial parameter~$\gamma_c=0$\,.

Conservation of fluid provides a third neutral mode in the dynamics. 
Thus, when~$\gamma=\tan\theta=\partial_x=\partial_y=\gamma_c=w_f=0$, a four parameter subspace of equilibria exists corresponding to  some uniform lateral shear, turbulent mean, flow, $(u,v)\propto z/h+c_u$, some turbulent mean concentration $c\propto 1-z/h$, on a fluid of any constant fluid depth~$h$. 
For large enough lateral length scales, these equilibria occur independently at each location~$x$ and~$y$~\cite[e.g.]{Roberts1988, Roberts:2008uq} and hence the subspace of equilibria are in effect parameterised by~$\uu(x,y)$, $\vv(x,y)$, $\cc(x,y)$ and~$h(x,y)$.

Provided we can treat lateral derivatives~$\partial_x$ and~$\partial_y$ as perturbing influences, that is provided solutions vary slowly enough in $x$~and~$y$, centre manifold theorems~\cite[e.g.]{Roberts1988, Chicone:2006fk} assure us of three vitally important properties:
\begin{itemize}
\item  this subspace of equilibria are perturbed to a slow manifold, whereon the evolution is slow, that exists for a finite range of gradients~$\partial_x$ and~$\partial_y$, and parameters~$\gamma$, $\gamma_c$, $\tan\theta$ and $w_f$, and which may be parameterised by the depth-averaged lateral velocities~$\uu(x,y,t)$ and~$\vv(x,y,t)$, the depth-averaged concentration~$\cc(x,y,t)$, and the local fluid depth~$h(x,y,t)$;
\item the slow manifold attracts solutions from all nearby initial conditions provided the spectrum of the linearised dynamics is suitable; 
\item and that a formal power series in the parameters~$\gamma$, $\gamma_c$, $\tan\theta$, $w_f$ and gradients~$\partial_x$ and~$\partial_y$ approximate the slow manifold to the same order of error as the order of the residuals of the governing differential equations.
\end{itemize}
That is, the theorems support the existence, emergence, and construction of slow manifold models such as~\eqref{Intro:model}.
\footnote{This general type of argument has recently been made rigorous in one lateral dimension by \cite{Roberts2013a}.}
This support occurs in a finite domain in parameter space, and we assume the finite domain is big enough to include interesting cases of the physically relevant \(\gamma=\gamma_c=1\) and finite slope~\(\tan\theta\) and falling velocity~\(w_f\).

\subsection{Linear dynamics of the system}
\label{smag:Tlinear}

The linear dynamics of the system support the application of centre manifold theory.
For the flat bed of~$b=\text{constant}$, and with~$\tan\theta=\gamma=\gamma_c=w_f=0$, the base problem~\eqref{sec:mass}--\eqref{sed:bcbedfG} has the equilibrium of a shear flow which is, in terms of the stretched vertical coordinate $Z=(z-b)/h$, 
\begin{align}&
h=\text{constant}\,,
\quad
u=2\uu\frac{Z+c_u}{1+2c_u}\,,
\quad
v=2\vv\frac{Z+c_u}{1+2c_u}\,,
\quad
w=0\,,\label{smag:linearSol1}
\\&
p=h(1-Z)\,,
\quad
c=2\cc(1-Z)\,,
\quad
\nu=c_th^2\ros=c_th\bq\frac{\sqrt{2}}{1+2c_u}\,,
\label{smag:linearSol}
\end{align}
where~$\uu$ and~$\vv$ are the depth-averaged lateral velocities, $\cc$ is the depth-averaged concentration, and~$\bq=\sqrt{\uu^2+\vv^2}$ is the mean speed.
Environmental turbulent flows have eddy viscosity variations.
In this linear analysis, we assume the eddy viscosity~$\nu$ is effectively constant.


Then we consider the dynamics of the~\pde{}s~\eqref{sec:mass}--\eqref{sec:con} linearised in the small perturbation fields $(h',u',v',w',c',p')$ about each of these equilibria.
Because environmental turbulent fluids have very large lateral scales compared with the depth the lateral variations are very slow.
As the lateral variations vary slowly they do not affect the dominant linear process.
Thus we also treat the lateral derivatives~$\partial_x$ and~$\partial_y$  as `small operators' in this linearisation.
Thus, the~\pde{}s~\eqref{sec:mass}--\eqref{sec:con} and the boundary conditions~\eqref{smag:nopenetrate}--\eqref{sed:bcbedfG}, in the linearisation that effectively $\partial_x=\partial_y=\tan\theta=\gamma=\gamma_c=w_f=0$, result in the linear problem
\begin{subequations}
\begin{align}&
\D{z}{w'}=0\,,
\label{smag:Lp1}\\&
\D{t}{u'}+w'\D zu=\nu\DD{z}{u'}\,
\label{smag:Lp2}\\&
\D{t}{v'}+w'\D zv=\nu\DD{z}{v'}\,,
\label{smag:Lp3}\\&
\D{t}{w'}=-\D{z}{p'}+\nu\DD{z}{w'}\,,
\label{smag:Lp4}\\&
\D{t}{c'}+w'\D zc=\nu\DD{z}{c'}\,,
\label{smag:Lp5}\\&
w'=0\quad\text{on}\quad Z=0\,,
\label{smag:Lp6}\\&
u'=c_uh'\D zu+c_uh\D{z}{u'}\quad\text{on}\quad Z=0\,,
\label{smag:Lp7}\\&
v'=c_uh\D{z}{v'}+c_uh'\D zv\quad\text{on}\quad Z=0\,,
\label{smag:Lp8}\\&
\D{t}{h'}=w'\quad\text{on}\quad Z=1\,,
\label{smag:Lp9}\\&
-p'+2\nu\D{z}{w'}=0\quad\text{on}\quad Z=1\,,
\label{smag:Lp10}\\&
\nu\D{z}{u'}=\frac{\sqrt{2}c_t\bq}{(1+c_u)(1+2c_u)}u'\quad\text{on}\quad Z=1\,,
\label{smag:Lp11}\\&
\nu\D{z}{\nu'}=\frac{\sqrt{2}c_t\bq}{(1+c_u)(1+2c_u)}v'\quad\text{on}\quad Z=1\,,
\label{smag:Lp12}\\&
c'=0\quad\text{on}\quad Z=1\,,
\label{smag:Lp13}\\&
\D{z}{c'}+\frac{c'}{h}=0\quad\text{on}\quad Z=0\,.
\label{smag:Lp14}
\end{align}
\label{smag:Lp}%
\end{subequations}
Equations~\eqref{smag:Lp1} and~\eqref{smag:Lp6} determine there is no vertical velocity, $w'=0$.
Equation~\eqref{smag:Lp9} implies the free surface perturbation $h'=\text{constant}$, which corresponds to the freedom already in~\eqref{smag:linearSol1} so without loss of generality we here set~$h'=0$.
Equations~\eqref{smag:Lp4} and~\eqref{smag:Lp10} implies no change to the hydrstatic pressure, $p'=0$.
\pde{}s~\eqref{smag:Lp2}, \eqref{smag:Lp3} and~\eqref{smag:Lp5}
indicate fields~$u'$, $v'$ and~$c'$ have solutions in the form 
\begin{equation}
(u',v',c')\propto[A\sin(kz)+B\cos(kz)]\exp(\lambda t)\,,
\label{sec:solF}
\end{equation}
where~$k$ is a nondimensional vertical wavenumber. 
Substitute these solution forms~\eqref{sec:solF} into \pde{}s~\eqref{smag:Lp2}, \eqref{smag:Lp3} and~\eqref{smag:Lp5}, and obtain $\lambda=-\nu k^2$.
Substituting the solution forms~\eqref{sec:solF} into boundary conditions~\eqref{smag:Lp11}--\eqref{smag:Lp14} leads to two separate conditions for the velocity fields and concentration fields, respectively
\begin{equation}
\tan k=\frac{k}{1+c_u(1+c_u)k^2}
\quad\text{and}\quad
\tan k=hk\,.\label{smag:kkk}
\end{equation}
The characteristic equations~\eqref{smag:kkk} have the non-zero wavenumbers $k>\pi$, which implies the leading non-zero eigenvalue $-\nu k^2<-\nu\pi^2$. 
In addition to these negative eigenvalues, there are zero eigenvalues  corresponding to the freedom to vary the fluid depth~$h$, depth-averaged velocities~$\uu$ and~$\vv$ and depth-averaged concentration~$\cc$.
Thus, there is a spectral gap between the eigenvalues $\lambda=0$ and $\lambda<-\nu\pi^2$.
Centre manifold theory~\cite[e.g.]{Roberts1988, Potzsche:2006uq} then supports the existence and emergence of a slow manifold of large lateral scale in the three dimensional turbulent fluid and sediment system.

\subsection{Reduced model of the flow and sediment dynamics}
\label{sec:model}

This section focusses on interpreting the slow manifold of the leading order model of the turbulent flow and suspended sediment.

Instead of depth-averaging equations, the centre manifold approximation theorem underlies an iterative construction of a slow manifold that resolves the turbulent and sediment interactions within the fluid layer. 
The computer algebra program listed in Appendix~\ref{caSmag} constructs the slow manifold of the turbulent flow and sediment system.
The program derives evolution equations for the water depth~$h(x,y,t)$, the depth-averaged lateral velocities $\uu(x,y,t)$~and~$\vv(x,y,t)$, and the depth-averaged concentration~$\cc(x,y,t)$. 

The order of error in the construction is phased in terms of the small parameters.
Here the small parameters are the lateral derivatives~$\partial_x$ and~$\partial_y$, the small mean slope~$\tan\theta$, the falling velocity~$w_f$, and the artificial parameters~$\gamma$ and~\(\gamma_c\).
Generally we report results to errors~$\mathcal O(\partial_x^{p/2}+\partial_y^{p/2}+\tan^{p/2}\theta+w_f^p+\gamma^p)$ for some prescribed exponent~$p$ (where, for example, a term with factor~\(w_f^m\gamma^n\) is~\Ord{w_f^p+\gamma^q} if \(m/p+n/q\geq1\)).
Further, the artificial small parameter~$\gamma_c$ is introduced to integrate the sediment dynamics into the theoretical support for a slow manifold.
In order to ensure the sediment dynamics are accurate we construct the slow manifold to higher orders in the parameter~$\gamma_c$.

\cite{Cao2014} [\S2.4.3] showed that coefficients of power series in~$\gamma$ converge quickly in the turbulent fluid flow systems.
Computations with the sediment equations indicate that the dependence upon~$\gamma_c$ also converges reasonably quickly.
Truncating to errors $\mathcal O(\partial_x^{3/2}+\partial_y^{3/2}+\tan^{3/2}\theta+\gamma^3+w_f^3,\gamma_c^5)$, the computer algebra in Appendix~\ref{caSmag} derives the evolution of the depth-averaged concentration~$\cc(x,y,t)$:
\begin{subequations}
\begin{align}
\D t\cc&={}\cdots
\nonumber\\&{}
-\left(\uu\D{x}{\cc}+\vv\D{y}{\cc}\right)\left(0.893+0.054\gamma_c+0.054\gamma_c^2+0.013\gamma_c^3-0.007\gamma_c^4\right)
\label{sed:cl4}\\&
{}+\frac{\cc}{h}\left(\uu\D{x}{h}+\vv\D{y}{h}\right)(0.052-0.024\gamma_c-0.027\gamma_c^2-0.007\gamma_c^3+0.003\gamma_c^4)
\label{sed:cl6}\\&
{}+\mathcal O(\partial_x^{3/2}+\partial_y^{3/2}+\tan^{3/2}\theta+\gamma^3+w_f^3,\gamma_c^5)\,.
\nonumber
\end{align}
\label{sed:cl}%
\end{subequations}
These power series converge quickly enough in~$\gamma_c^5$ to reasonably evaluate at $\gamma_c=1$ to give accurate coefficients in the evolution equations.

Executing the computer algebra in Appendix~\ref{caSmag} and evaluating at~$\gamma=\gamma_c=1$ leads to the following evolution equations for the turbulent flow and sediment system in terms of the depth~$h(x,y,t)$, depth-averaged lateral velocities $\uu(x,y,t)$ and $\vv(x,y,t)$ and depth-averaged concentration~$\cc(x,y,t)$.
The equations here are complicated due to the methodology systematically resolving all the intricate microscale physical interactions.
\begin{subequations}
\begin{align}
\D{t}{h}&\approx{}-\left(\D{x}{h\bar u}+\D{y}{h\bar v}\right)\,,
\label{smag:h}\\
\D{t}{\bar u}&\approx{}-0.00293\frac{\bar u\bq}{h}+0.993\left[\tan\theta-\D{x}{(h+b)}\right]
\nonumber\\&
{}-1.025\uu\D x\uu-1.017\vv\D y\uu-0.00817\left(\frac{\uu^2}{h}\D xh-\frac{\uu\vv}{h}\D yh\right)
\nonumber\\&
{}+0.0941\frac{\bq}{h}\left[\D{x}{}\left(h^2\D{x}{\uu}\right)+\D{y}{}\left(h^2\D{y}{\uu}\right)\right]
\nonumber\\&
{}+0.0839\frac{\uu^2-\vv^2}{h\bq}\left[\D{x}{}\left(h^2\D{x}{\uu}\right)-\D{y}{}\left(h^2\D{y}{\uu}\right)\right]
\nonumber\\&
{}+0.00257(s-1)\frac{\uu\cc\bq}{h}-0.298(s-1)h\D x\cc\,,
\label{smag:u}\\
\D{t}{\vv}&\approx{}-0.00293 \frac{\vv\bq}{h}-0.993\D{y}{(h+b)}
\nonumber\\&
{}-1.025\vv\D y\vv-1.017\uu\D x\vv-0.00809\left(\frac{\uu\vv}{h}\D xh-\frac{\vv^2}{h}\D yh\right)
\nonumber\\&
{}+0.0941\frac{\bq}{h}\left[\D{x}{}\left(h^2\D{x}{\vv}\right)+\D{y}{}\left(h^2\D{y}{\vv}\right)\right]
\nonumber\\&
{}+0.0839\frac{\uu^2-\vv^2}{h\bq}\left[\D{x}{}\left(h^2\D{x}{\vv}\right)-\D{y}{}\left(h^2\D{y}{\vv}\right)\right]
\nonumber\\&
{}+0.00257(s-1)\frac{\vv\cc\bq}{h}-0.298(s-1)h\D y\cc\,,
\label{smag:v}\\
\D{t}{\cc}&\approx
{}-\frac{w_f}{h}\left(0.938+28.9\frac{w_f}{\bq}\right)\cc+\frac{w_f}{h}\left(0.984-51.3\frac{w_f}{\bq}\right)c_{ae}
\nonumber\\&
{}-\left(1.01-3.09\frac{w_f}{\bq}\right)\left(\uu\D x\cc+\vv\D y\cc\right)
\nonumber\\&
{}+0.0331\frac{\bq}{h}\left[\D{x}{}\left(h^2\D{x}{\cc}\right)+\D{y}{}\left(h^2\D{y}{\cc}\right)\right]
\nonumber\\&
{}+0.0271\frac{\uu^2-\vv^2}{h\bq}\left[\D{x}{}\left(h^2\D{x}{\cc}\right)-\D{y}{}\left(h^2\D{y}{\cc}\right)\right]
\label{sed:cH}
\end{align}
\end{subequations}%
Equation~\eqref{smag:h} is a direct consequence of conservation of fluid.
The momentum equations~\eqref{smag:u}--\eqref{smag:v} include the effects of drag $\uu\bq/h$, advection, such as~$\uu\D{x}{\uu}$, turbulent dissipation, gravitational forcing $\tan\theta-\D{x}{(h+b)}$, and pressure gradients established by the suspended sediment~$h\D{x}{\cc}$.
The sediment concentration equation~\eqref{sed:cH} contains the equilibration of vertical sediment distribution due to terms such as~$w_f\cc/h$,  advection such as~$\uu\D{x}{\cc}$, and turbulent dispersion effects. 

Although equations~\eqref{smag:h}--\eqref{sed:cH} are expressed in terms of depth-averaged lateral velocities and depth-averaged concentration, they are derived not by depth-averaging, but instead by systematically accounting for interaction between vertical profiles of the velocity and concentration, the stress, bed drag, lateral space variations and bed topography. 
The form of coefficients in equations~\eqref{smag:h}--\eqref{sed:cH} are supported by dynamical systems theory: the detail in the equations reflects that a slow manifold is in principle composed of exact solutions of the Smagorinski dynamics~\eqref{sec:mom} and convection-diffusion equation~\eqref{sec:con}, and hence accounts for all interactions up to a given order of analysis no matter how small the numerical coefficient in the interactions.

The sediment model~\eqref{sed:cH} includes all the dominated terms in established modelling~\cite[e.g.]{Wu:2004uq,Duan2004,Duan:2006fk}.
For example, \cite{Duan2004} derived the following depth-averaged advection-diffusion equation of suspended sediment:
\begin{align}
\D{t}{\cc}&={}-\frac{w_f}{h}(\cc-c_{ae})-\uu\D{x}{\cc}-\vv\D{y}{\cc}+0.13h\bq\left(\DD{x}{\cc}+\DD{y}{\cc}\right)\,,
\label{sed:establish}
\end{align}
which consists of the effects of vertical distribution~$w_f/h(\cc-c_{ae})$, advection $\uu\D{x}{\cc}$ and~$\vv\D{y}{\cc}$, and dispersion~$h\bq\DD{x}{\cc}$ and~$h\bq\DD{y}{\cc}$.
But the derived model~\eqref{sed:cH} also includes more subtle effects, which could be important for suspended sediment in complex flow regimes.
The model~\eqref{sed:cH} further gives modifications in presence of the ratio~$w_f/\bq$ due to different distributions of sediment in the vertical for the different levels of turbulent mixing characterised by the mean flow speed~$\bq$.
The coefficients in equation~\eqref{sed:cH} are a little different to the established model~\eqref{sed:establish}.
Take the advection term~$\uu\D{x}{\cc}$ for example, the established model~\eqref{sed:establish} has the coefficient~$1$, whereas our coefficient of such terms is $1.01-3.09w_f/\bq$.
Physically, a higher falling velocity~$w_f$ means sediment concentrates more near bed where the mean advection velocity is lower, and hence net transport will be slower.
Our model~\eqref{sed:cH} has smaller dispersion coefficent compared with the model~\eqref{sed:establish}.

\section{Cross-sectional structures of the flow and sediment}
\label{sec:Ver}

The centre manifold approach does not impose a specific cross-sectional velocity distribution as done by other methods, instead our approach empowers the sediment and Smagorinski turbulent equations to determine the cross-sectional structures in out-of-equilibrium dynamics.
Recall that the locally stretched vertical coordinate $Z=(z-b)/h$\,.
This section concentrates on the vertical distribution of the lateral velocity~$u(Z)$ and concentration~$c(Z)$ in steady flow, and compare our prediction of the lateral velocity~$u(Z)$ with analogous published experimental data~\cite[e.g.]{Schultz:2007uq, Schultz:2013fk, Celik:1988fk}.

\subsection{Distribution of the suspended sediment}
\label{sec:VerC}

We concentrate on the vertical distribution of the suspended sediment~$c(Z)$ in the slow manifold of the system.
Computer algebra in Appendix~\ref{caSmag} derives the physical field of sediment concentration~$c$ in terms of the depth~$h(x,y,t)$, depth-averaged velocities~$\uu(x,y,t)$ and~$\vv(x,y,t)$, depth-averaged concentration~$\cc(x,y,t)$ and stretched local normal coordinate $Z=(z-b)/h$ on the slow manifold, evaluating at $\gamma=\gamma_c=1$:
\begin{subequations}
\begin{align}
c(Z)&={}\cc\left(0.985+0.0422Z-0.00756Z^2-0.0139Z^3\right)
\label{sed:conV1}\\&
{}+\cc\frac{w_f}{\bq}\left(28.36-5.156Z-77.34Z^2\right)
\label{sed:conV2}\\&
{}+\cc\frac{w_f^2}{\bq}\left(-0.430-0.430Z+2.578Z^2-0.859Z^3\right)
\label{sed:conV3}\\&
{}+c_{ae}\frac{w_f}{\bq}\left(56.72-166.3Z+77.34Z^2+2.578Z^3\right)
\label{sed:conV4}\\&
{}+c_{ae}\frac{w_f^2}{\bq}\left(-0.430+1.074Z-0.430Z^3\right)
\label{sed:conV5}\\&
{}+\frac{h}{\bq}\left(\uu\D x\cc+\vv\D y\cc\right)\left(2.578+0.921Z-17.68Z^2+11.42Z^3\right)
\label{sed:conV6}\\&
{}+\frac{h\cc}{\bq}\D x\uu\left(-0.17+0.449Z-0.0392Z^2-1.322Z^3\right)
\label{sed:conV7}\\&
{}+\frac{h\cc}{\bq}\D y\vv\left(-1.774+1.244Z+2.471Z^2-1.486Z^3\right)
\label{sed:conV8}\\&
{}+\frac{\cc}{\bq}\left(\uu\D xh+\vv\D yh\right)\left(-0.17+0.449Z-0.0392Z^2-1.322Z^3\right)
\label{sed:conV9}\\&
{}+\frac{\cc}{\bq^3}\left(\uu\D xb+\vv\D yb\right)\left(0.918-0.809Z_2.549Z^2+1.343Z^3\right)
\label{sed:conV10}\\&
{}+\frac{\cc}{\bq^3}\left(\uu\D xb+\vv\D yb\right)\left(0.918-0.809Z_2.549Z^2+1.343Z^3\right)
\label{sed:conV11}\\&
{}-\tan\theta\,\frac{h\uu\cc}{\bq^3}\left(0.918-0.809Z+2.549Z^2+1.343Z^3\right)
\label{sed:conV12}\\&
{}+\Ord{\partial_x^{3/2}+\partial_y^{3/2}+\tan^{3/2}\theta+w_f^3+\gamma^3,\gamma_c^5}\,.
\nonumber
\end{align}
\label{sed:conV}%
\end{subequations}
The vertical sediment distribution~\eqref{sed:conV} describes the low-order shape of the slow manifold in state space.
Physically this equation describes the details of the suspended sediment concentration in out-of-equilibrium flow. 
The terms in equation~\eqref{sed:conV} have physical interpretations.
For example, the line~\eqref{sed:conV1} approximates the mean concentration, together with the lines~\eqref{sed:conV2}--\eqref{sed:conV5}, which forms the basic distribution of the concentration~$c(Z)$ in the vertical in the presence of the depth-averaged concentration~$\cc$, the equilibrium bed concentration~$c_{ae}$, the falling velocity~$w_f$, and the mean fluid speed~$\bq$.
The lines~\eqref{sed:conV6}--\eqref{sed:conV9} describe the effect by advection on the vertical distribution of the concentration.
Lines~\eqref{sed:conV10}--\eqref{sed:conV11} describe modifications due to the change of the bed topography.
The line~\eqref{sed:conV12} describes modifications due to lateral flow down a slope in its affect on the vertical distribution of the suspended sediment.

Now we explore the distribution of the suspended sediment~$c(Z)$ in steady flows.
Consider the turbulent flow of constant depth~$H=1$ with suspended sediment flowing on a flat mean bed of constant slope~$\tan\theta$; that is, the mean bed~$b=0$.
We consider the concentration fraction is small, $\cc<0.01$.
In this regime, in \eqref{smag:u}--\eqref{smag:v} the terms $\uu\cc\bq/h$ and $\vv\cc\bq/h$ are negligible.
Then the evolution equations~\eqref{smag:h}--\eqref{smag:v} predicts the equilibrium $U=18.7\tan^{1/2}\theta$ and $V=0$, so the mean speed $\bq=U=18.7\tan^{1/2}\theta$.

The nondimensional equilibrium reference concentration~$c_{ae}$ on the mean bed $z=b$ in steady flow is determined by the particle size and the slope.
We approximate the nondimensional shear velocity~$u_f=\bq/C'$, where~$C'=18\log(4/d)$ is the Chezy coefficient~\cite[e.g.]{Rijn:1984fk}.
Thus, equation~\eqref{sed:cae1} gives the nondimensional equilibrium reference concentration
\begin{equation}
c_{ae}=3.26\tan^{1.5}\theta\,\frac{1}{d^{0.8}(1.39-\log d)^3}\,,
\label{sed:cae}
\end{equation}
which only depends on the nondimensional mean slope~$\tan\theta$ and the nondimensional mean particle size~$d$.

For the steady flow, the concentration $c(Z)$ in equation~\eqref{sed:conV} then reduces to
\begin{align}
c(Z)&={}\cc\left(0.985+0.0422Z-0.00756Z^2-0.0139Z^3\right)
\nonumber\\&
{}+\cc\frac{w_f}{\bq}\left(28.36-5.156Z-77.34Z^2\right)
\nonumber\\&
{}+\cc\frac{w_f^2}{\bq}\left(-0.430-0.430Z+2.578Z^2-0.859Z^3\right)
\nonumber\\&
{}+c_{ae}\frac{w_f}{\bq}\left(56.72-166.3Z+77.34Z^2+2.578Z^3\right)
\nonumber\\&
{}+c_{ae}\frac{w_f^2}{\bq}\left(-0.430+1.074Z-0.430Z^3\right).
\label{sed:conVS}
\end{align}
Equation~\eqref{sed:conVS} shows that the concentration~$c(Z)$ depends on the vertical coordinate~$Z$, the falling velocity~$w_f$, the depth-averaged concentration~$\cc$, the equilibrium reference concentration~$c_{ae}$ and the mean flow speed~$\bq$.
The falling velocity~$w_f$ and the equilibrium reference concentration~$c_{ae}$ vary with the mean particle size~$d$ and the mean slope~$\tan\theta$ according to equation~\eqref{sed:wf} and \eqref{sed:cae}, respectively.
Thus, the concentration profile~$c(Z)$ depends on the vertical coordinate~$Z$, the mean particle size~$d$ and the mean slope~$\tan\theta$.

\begin{figure}
\centering
\begin{tabular}{c@{}c}
\rotatebox{90}{\hspace{20ex}$Z$} &
\includegraphics{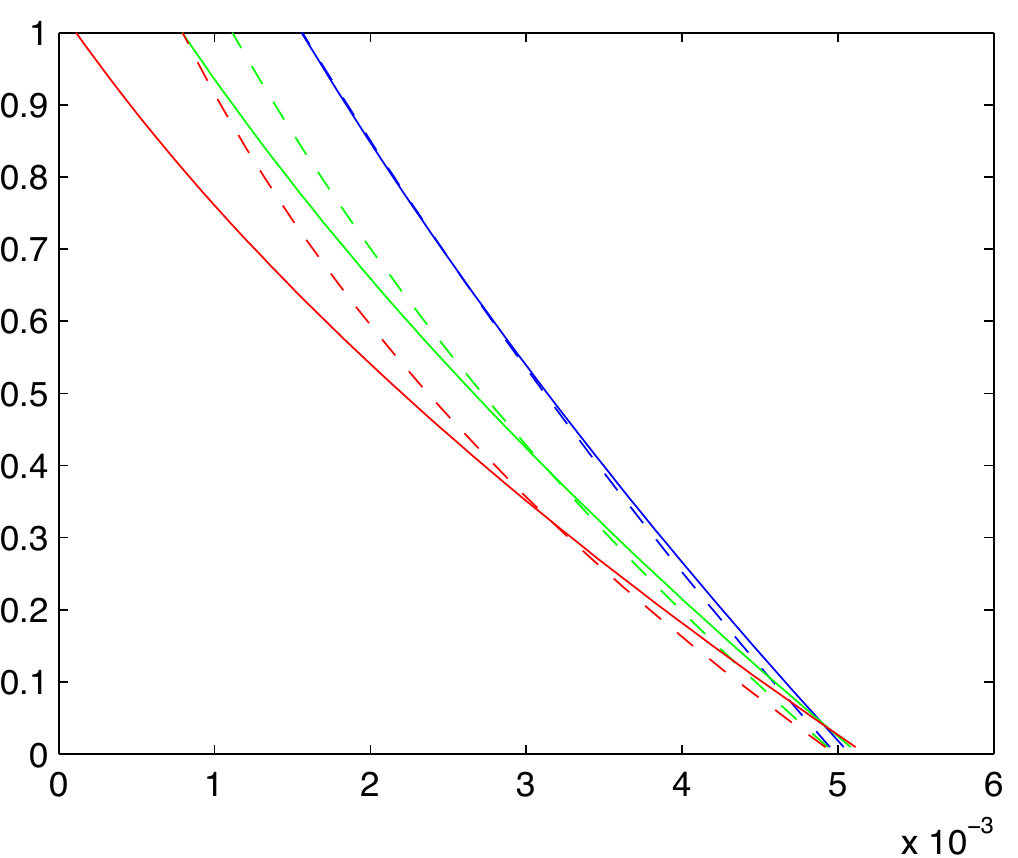}\\
& $c(Z)$
\\&
$\color{blue}\line(1,0){20}$~$d=0.6\times10^{-4}$ \ 
$\color{green}\line(1,0){20}$~$d=1\times10^{-4}$ \ 
$\color{red}\line(1,0){20}$~$d=1.5\times10^{-4}$ \ 
\end{tabular}
\caption{Profiles (line curves) of the suspended sediment concentration~$c(Z)$ in the vertical for three different nondimensional mean particle size~$d$ according to equation~\eqref{sed:conVS}.
The mean slope $\tan\theta=0.01$.
The equilibrium reference concentration~$c_{ae}=0.005$.
The dash curves are the corresponding steady analytical approximation~\eqref{sed:CZS}.}
\label{sedCZ}
\end{figure}

\begin{figure}
\centering
\begin{tabular}{c@{\ }c}
\rotatebox{90}{\hspace{20ex}$Z$} &
\includegraphics{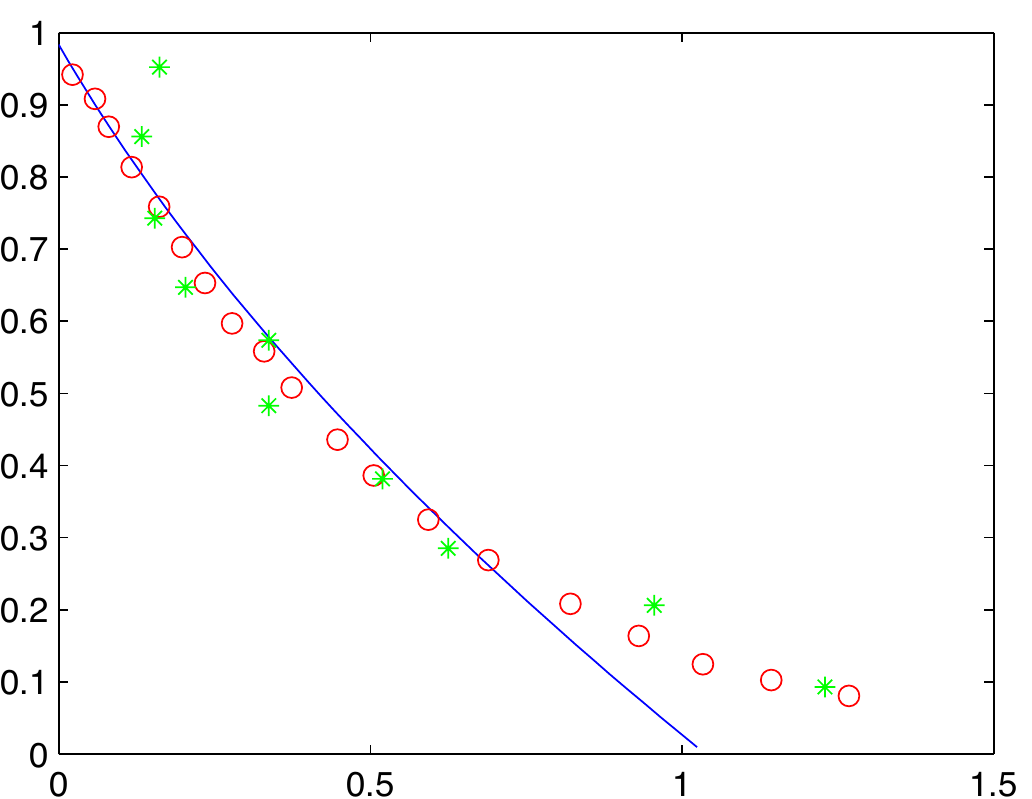}\\
& $\cc/c_{ae}$
\end{tabular}
\caption{Vertical distribution of the suspended sediment: (blue curve) from equation~\eqref{sed:conVS}; (red circle) the numerical prediction by~\cite{Celik:1988fk}; and (green stars) the corresponding experimental data used by~\cite{Celik:1988fk}.
The nondimensional mean particle size~$d=1.65\times10^{-4}$ and then the falling velocity~$w_f=0.0161$ in our simulation.}
\label{sedCZcom}
\end{figure}

Figure~\ref{sedCZ} depicts the profiles of the nondimensional suspended sediment concentration~$c(Z)$ in the vertical for four different nondimensional mean particle size~$d$.
The bed has a mean slope $\tan\theta=0.01$.
The nondimensional equilibrium reference concentration~$c_{ae}\approx0.005$ for the mean particle size $d>6\times10^{-5}$.
For small nondimensional mean particle size~$d$, the nondimensional concentration~$c(Z)$ is approximately linear in the vertical coordinate~$Z$.
In this steady flow, the approximation of the eddy diffusivity is
\begin{align}
\epsilon_s(Z)&\approx{} \bq(0.00628-0.00269Z-0.000733Z^2)
\nonumber\\&
{}+\tan\theta\,\frac{1}{\bq}(0.00978-0.2605Z+0.247Z^2)\,.
\label{sed:conNu}
\end{align}
For an indicative comparison, we integrate equation~\eqref{sec:con} from the bottom to the free surface in the steady flow, and obtain an approximation
\begin{equation}
c(Z)\approx c_{ae}\left[\frac{5.29+Z}{3.26(1.62-Z)}\right]^{-197.45w_f/\bq}\,.
\label{sed:CZS}
\end{equation}
The dash curves in Figure~\ref{sedCZ} depict the approximation~\eqref{sed:CZS}.
When the nondimensional mean particle size~$d$ is small, the computed vertical distribution~\eqref{sed:conVS} is approximately the same with the approximation~\eqref{sed:CZS}.
When the nondimensional mean particle size~$d$ increases, there is a difference between the computed vertical distribution~\eqref{sed:conVS} and the approximation~\eqref{sed:CZS} at the upper flow.

Figure~\ref{sedCZcom} plots the vertical distribution of the suspended sediment from equation~\eqref{sed:conVS}~(blue curve), the numerical prediction~(red circles) by~\cite{Celik:1988fk}, and the corresponding experimental data~(green stars) used by~\cite{Celik:1988fk}.
\cite{Celik:1988fk} calculated the suspended sediment transport in unidirectional channel flow, where they used the nondimensional variables of fluid depth~$H=1$, mean velocity $U\approx1.8$, mean particle size $d=1.65\times10^{-4}$ and falling velocity $w_f=0.0165$. 
Our simulation~(blue curve) agrees with the numerical prediction by~\cite{Celik:1988fk} except at the bottom.
This difference at the bottom is because we have small entrainment at the mean bed.
Our simulation~(blue curve) is good enough to predict the experimental data~(green stars).
The trends of the suspended sediment concentration qualitatively agrees with other published experimental measurements~\cite[e.g.]{Cellino:1999fk,Yoon:2005fk}.

\subsection{Vertical distribution of the velocity in steady flow}
\label{sec:VerV}

This section reports on the vertical distribution of the lateral velocity~$u(Z)$ in the flow and sediment system, \(v(Z)\)~is similar.
Computer algebra in Appendix~\ref{caSmag} derives the following physical flow field of lateral velocity~$u$ in term of stretched local normal coordinate $Z=(z-b)/h$ on the slow manifold, evaluated at the physical $\gamma=\gamma_c=1$:
\begin{align}
u(Z)&={}\uu(0.816+0.445Z-0.0916Z^2-0.0307Z^3-0.00383Z^4-0.000418Z^5)
\nonumber\\&
{}+\tan\theta\,\frac{h}{\bq}(2.208+1.204Z-14.31Z^2+8.069Z^3-1.569Z^4
\nonumber\\&\quad{}
+0.954Z^5+0.586Z^6+0.119Z^7)
\nonumber\\&
{}+\frac{h\uu}{\bq}\D x\uu(2.326+1.269Z-13.52Z^2+4.585Z^3+0.894Z^4
\nonumber\\&\quad{}
+0.783Z^5+0.533Z^6+0.118Z^7+0.0106Z^8)
\nonumber\\&
{}+\frac{h\vv}{\bq}\D y\uu(2.352+1.283Z-13.25Z^2+3.622Z^3+1.53Z^4
\nonumber\\&\quad{}
+0.708Z^5+0.543Z^6+0.129Z^7+0.0127Z^8)
\nonumber\\&
{}+\frac{h}{\bq}\left(\D xh+\D xb\right)(-2.208-1.204Z+14.31Z^2-8.069Z^3
\nonumber\\&\quad{}
+1.569Z^4-0.954Z^5-0.586Z^6-0.119Z^7)\,,
\nonumber\\&
{}+(s-1)\cc\uu\left(0.0167+0.009Z-0.107Z^2+0.0439Z^3+0.0173Z^4\right)
\nonumber\\&
{}+\Ord{\partial_x^{3/2}+\partial_y^{3/2}+\tan^{3/2}\theta+w_f^3+\gamma^3,\gamma_c^5}\,.
\label{smag:uZ}
\end{align}
Physically, equation~\eqref{smag:uZ} describes the vertical details of the lateral velocity~$u(Z)$ in terms of the vertical position~$Z$, fluid depth~$h$, depth-averaged velocities~$\uu$ and~$\vv$, depth-averaged concentration~$\cc$, and the slope $\tan\theta$ of the mean bed~$b$.

\cite{Schultz:2013fk} experimentally studied the smooth-wall turbulent channel flow with the Reynolds number~$\re$ up to~$300\,000$, and showed that the mean flow is approximately independent of the Reynolds number.
In earlier work, \cite{Schultz:2007uq} showed when the relative roughness, the ratio of the roughness height and the boundary-layer thickness, is small, the mean velocity shape for the rough and smooth walls are similar in the outer layer.

However, compared with large roughness in the environmental flows under consideration, the roughnesses in the experiments of \cite{Schultz:2007uq} are small---the ratio between the roughness height and the fluid depth was${}\approx1.5\times10^{-3}$.
To compare with the experimental data~\cite[e.g.]{Schultz:2007uq, Schultz:2013fk}, we derive the equilibria profiles and evaluate the lateral velocity at these equilibria.

\begin{figure}
\centering
\begin{tabular}{c@{\ }c}
\rotatebox{90}{\hspace{20ex}$Z$} &
\includegraphics{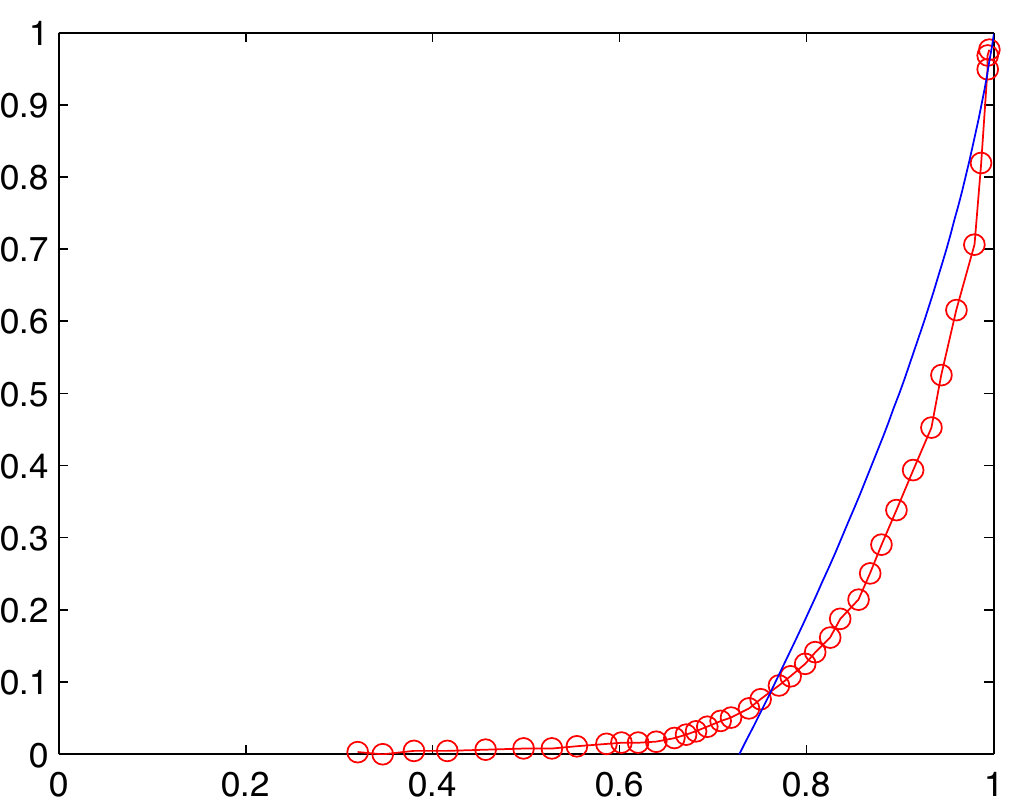}\\
& $u(Z)/u(1)$
\end{tabular}
\caption{Vertical distribution of the lateral velocity from: (blue curve) the approximation~\eqref{smag:uZS}; and (red circle curve) the experimental measurements~\cite[Fig.~2]{Schultz:2013fk}.
The ratio~$u(Z)/u(1)$ is independent of the slope~$\tan\theta$, where~$u(1)$ is the lateral velocity at the level~$Z=1$.
}
\label{smagUZC}
\end{figure}

Consider the turbulent flow with suspended sediment flowing on a flat mean bed of constant slope~$\tan\theta$; that is, the mean bed~$b=0$.
Consider the suspended sediment in steady flow of depth~$H=1$.
We consider the concentration fraction is small, $c_{ae},\cc<0.01$\,,
so that in \eqref{smag:u}--\eqref{smag:v} the terms $\uu\cc\bq/h$ and $\vv\cc\bq/h$ are negligible.
Then the evolution equations~\eqref{smag:h}--\eqref{smag:v} predicts the equilibrium $U=18.7\tan^{1/2}\theta$ and $V=0$, so the mean speed $\bq=U=18.7\tan^{1/2}\theta$.
For this steady flow, the lateral velocity~\eqref{smag:uZ} reduces to 
\begin{align}
\frac{u(Z)}{\sqrt{\tan\theta/c_t}}&\approx{}2.18+1.19Z-0.297Z^2-0.0533Z^3-0.0173Z^4
\nonumber\\&
{}-0.00366Z^5-0.00115Z^6-0.000089Z^7\,.
\label{smag:uZS}
\end{align}
Figure~\ref{smagUZC} compares the profile of the lateral velocity from the approximation~\eqref{smag:uZS}~(blue line curve) with the experimental measurements~(red circle curve) by~\cite{Schultz:2013fk}.
From the approximation~\eqref{smag:uZS}, the velocity ratio $u(Z)/u(1)$ is independent of the slope~$\tan\theta$.
Our prediction of the lateral velocity in the vertical agrees reasonably the experimental lateral velocity by~\cite{Schultz:2013fk}, except near the mean bed.
In environmental flows we expect that the large roughness of stones, roots and debris to typically break up any log boundary layer.
Thus, we do not resolve the turbulent log layer, and are interested in the dynamics determined by the relatively large scale of the fluid depth.

\begin{figure}
\centering
\begin{tabular}{c@{\ }c}
\rotatebox{90}{\hspace{15ex}$\tau_{xz}(Z)/\tan\theta$} &
\includegraphics{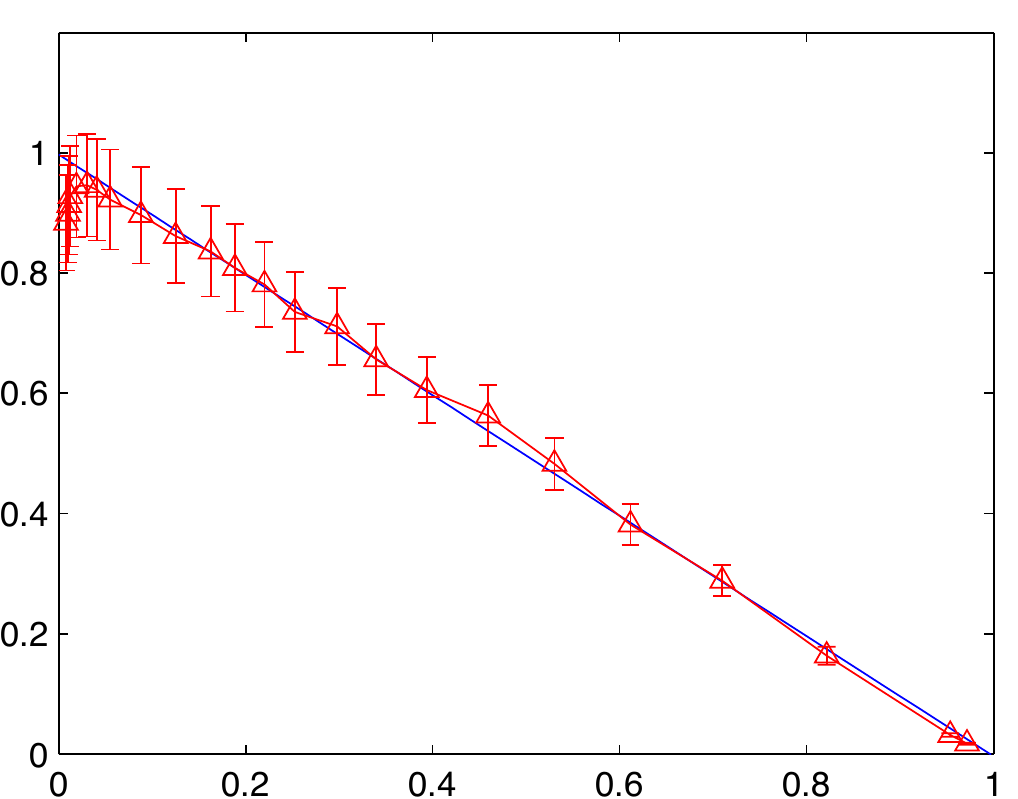}\\
& $Z$
\end{tabular}
\caption{Shear stress profile of: (blue curve) our approximation $\tau_z$; and (red triangle curve) the experimental measurements of flow over rough bed~\cite[Fig.~9]{Schultz:2007uq}.
The error bars show the~$\pm9\%$ uncertainty in their experiments. 
}
\label{smagSS}
\end{figure}

Shear stress arises in the turbulent fluid.
In this steady flow, we predict the shear stress~$\tau_{xz}$ has the near linear profile of
\begin{align}
\frac{\tau_{xz}(Z)}{\tan\theta}&\approx0.997-0.999Z+0.000284Z^2-0.00995Z^3
\nonumber\\&
{}+0.00776Z^4+0.000791Z^5+0.000072Z^6\,.\label{smag:nuZS}
\end{align}
Figure~\ref{smagSS} shows that our approximation of shear stress~$\tau_{xz}$~(blue curve)  approaches the experimental measurements~(red triangle curve) by~\cite{Schultz:2007uq}.
The shear stress is approximately straight due to the need to dissipate the near constant forcing of gravity. 
A difference occurs near the bed, because we do not resolve the log boundary layer.

\section{Simulating suspended sediment in large waves over a rippled bed}
\label{sec:Num}

This section explores the depth-averaged concentration in waves on an inclined rippled bed by the suspended sediment model~\eqref{sed:cH}, coupled with the modified momentum equations \eqref{smag:h}--\eqref{smag:v}.
Numerical results are qualitatively compared with the experimental measurements of the suspended sediment~\cite[e.g.]{Zedler:2006fk,KosYan:2007fk}.

Consider the fluid of depth~$h(x,y,t)$ with suspended sediment of depth-averaged concentration~$\cc(x,y,t)$ flowing down a rippled bed.
The fluid has the depth-averaged lateral velocities~$\uu(x,y,t)$ and~$\vv(x,y,t)$ along the bed.
Let the bed have nondimensional length~$L_x=100$ and width~$L_y=10$.
The mean slope of the bed is $\tan\theta=0.01$.
The green curve in Figure~\ref{sedPRhu} shows the ripples on the bed.
The ripples have the maximum nondimensional height~$0.4$ and length~$20$.
The bed has zero mean.

In simulations, all the variables are nondimenisonal.
For this pilot study, we consider the depth-averaged velocity~$\vv=0$ throughout.
The flow has mean equilibrium depth~$H=1$, then the mean equilibrium depth-averaged lateral velocity $U\approx1.86$ and mean equilibrium depth-averaged concentration $\cc\approx0.0035$.
The simulations initially impose a small perturbation $0.2\sin(2\pi/L_xx)$ to the equilibrium.
Recall the Froude number $U/\sqrt{gH}=1.86>1$ here.
Thus, we predict supercritical flow arises in the simulation.

\begin{figure}
\centering
\begin{tabular}{c@{}c}
\rotatebox{90}{\hspace{15ex} velocity~$\uu$} &
\includegraphics{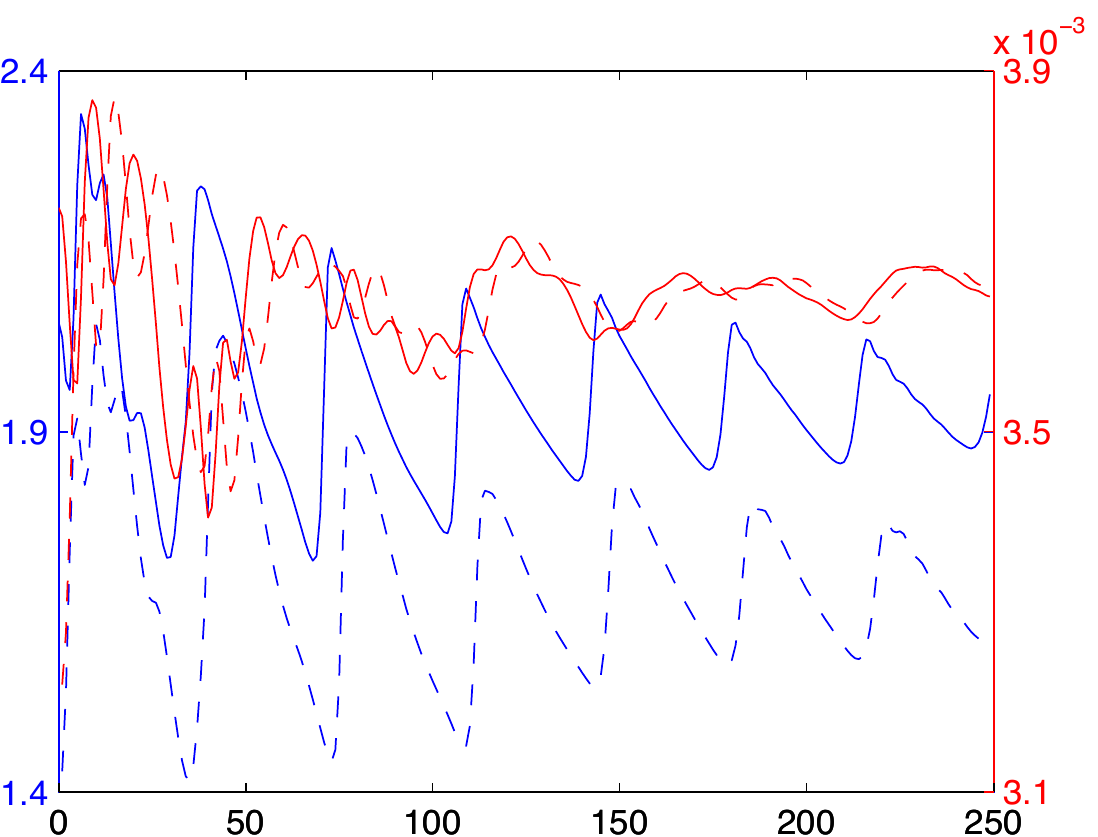}%
\rotatebox{90}{\hspace{15ex} concentration~$\cc$}\\
& $t$
\end{tabular}
\caption{Time series of the depth-averaged lateral velocity~$\uu$~(blue curves) and the depth-averaged concentration~$\cc$~(red curves) at the trough $x=50$~(line curves) and at the crest $x=60$~(dash curves) in Figure~\ref{sedPRhu}.
The mean particle size~$d=6\times10^{-5}$, so the falling velocity~$w_f=0.0097$ and the equilibrium reference concentration $c_{ae}=0.0057$.}
\label{sedPRt}
\end{figure}

Simulate the suspended sediment in the fluid flowing over the rippled bed of Figure~\ref{sedPRhu} by the developed model~\eqref{smag:h}--\eqref{sed:cH} with periodic boundary conditions in both $x$~and~$y$ directions for both the flow and bed.
Figure~\ref{sedPRt} plots the time series of the depth-averaged lateral velocity~$\uu$~(blue curves) and the depth-averaged concentration~$\cc$~(red curves) at a trough $x=50$~(line curves) and at a crest $x=60$~(dash curves).
The periodic depth-averaged velocity~$\uu$ indicates that large roll waves are generated on the free surface~\cite[e.g.]{Balmforth2004}.
The depth-averaged velocity~$\uu$ is bigger at the trough $x=50$~(blue line curve) than at the crest $x=60$~(blue dash curve), which indicates enhanced turbulent mixing arises at the trough~\cite[e.g.]{Zedler:2001fk}.
Then the turbulent mixing produces slightly bigger depth-averaged concentration~$\cc$ at the trough $x=50$~(red line curve) than at the crest $x=60$~(red dash curve).
\cite{Zedler:2006fk}, in their calculation of suspended sediment over rippled beds, found that far from the bed the concentration at the crest and trough are approximately~$\pi$ out of phase.
They commented that this phase lag is due to the vortex produced by the ripple near the trough.
\cite{KosYan:2007fk} experimentally and mathematically plotted the time series of the horizontal velocity and concentration under waves over sandy bottom, which showed that the concentration slightly lags the horizontal velocity to reach maximum. 
However, in our simulation,  no significant lag happens which appears due to our ripple not producing a vortex near the trough.

\begin{figure}
\centering
\begin{tabular}{c@{}c}
\rotatebox{90}{\hspace{20ex}$h\ \&\ \uu$} &
\includegraphics{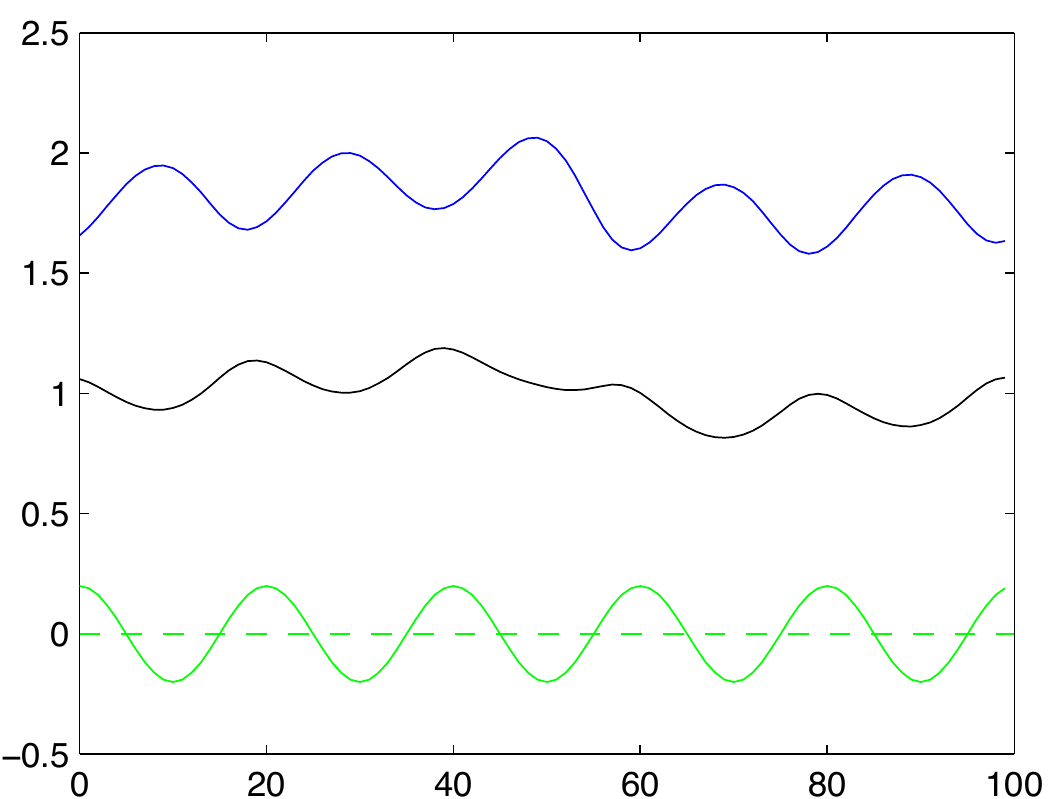}\\
& $x$
\end{tabular}
\caption{Plots of the depth~$h$~(black) and depth-averaged velocity~$\uu$~(blue) of the fluid flowing over the rippled bed~(green) at time $t=180$ in Figure~\ref{sedPRt}.
The dash line represents the zero mean bed.
This figure shows supercritical flow arises.}
\label{sedPRhu}
\end{figure}

\begin{figure}
\centering
\begin{tabular}{c@{}c}
\rotatebox{90}{\hspace{15ex}depth~$h$} &
\includegraphics{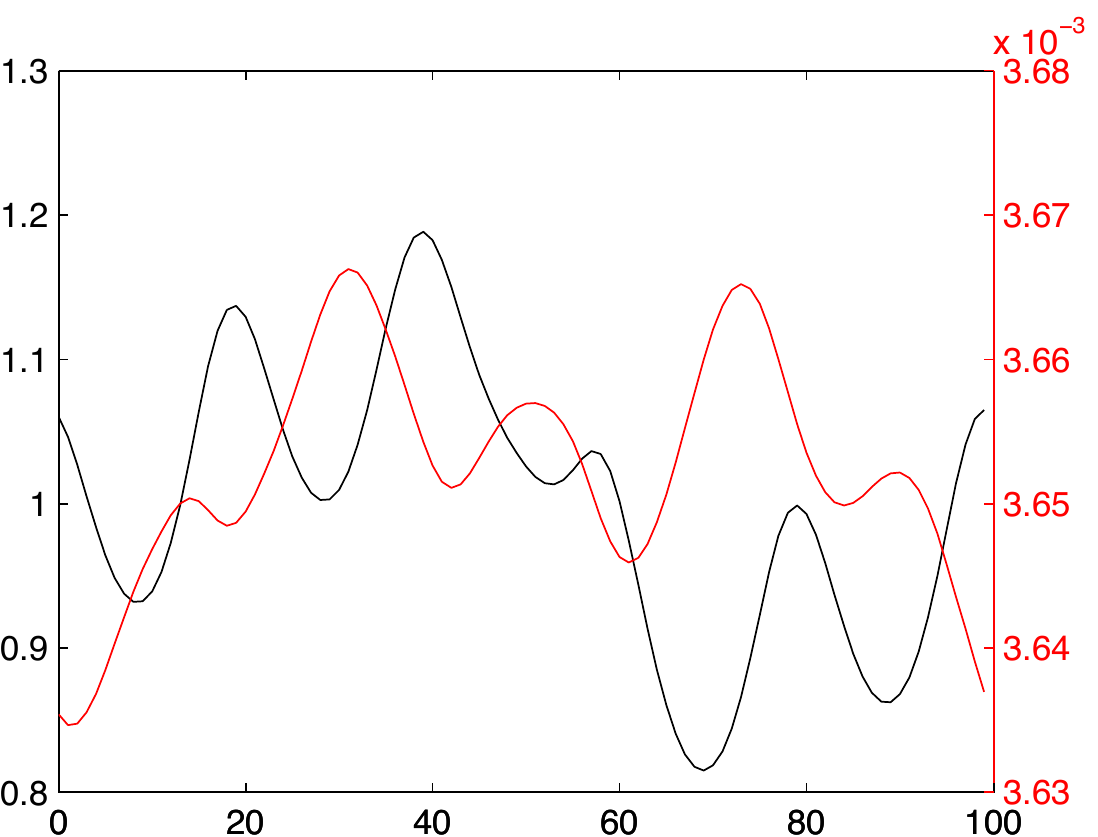}%
\rotatebox{90}{\hspace{15ex}\color{red}concentration~$\cc$}\\
& $x$
\end{tabular}
\caption{Plots of the depth~$h$~(black) and the depth-averaged concentration~$\cc$~(red) in the $x$~direction at time~$t=180$ in Figure~\ref{sedPRt}.
The depth-averaged concentration~$\cc$ is ahead to reach maximum over a ripple.}
\label{sedPRx}
\end{figure}

Figure~\ref{sedPRhu} plots the water depth~$h$~(black) and depth-averaged velocity~$\uu$~(blue) of the fluid flowing over the rippled bed~(green) in the $x$~direction at time $t=180$ in Figure~\ref{sedPRt}.
The dash green line represents the zero mean bed level.
Figure~\ref{sedPRhu} exhibits the supercritical flow as the fluid flowing over each ripple on the bed.
The depth~$h$ rises at the crest and the depth-averaged velocity~$\uu$ declines at the crest, which corresponds to the depth-averaged velocity~$\uu$ reaches minimum at the crest in Figure~\ref{sedPRt}.
Figure~\ref{sedPRx} plots the depth~$h$ and the depth-averaged concentration~$\cc$ in $x$-direction~at time $t=180$ in Figure~\ref{sedPRt}.
The depth-averaged concentration~$\cc$ is approximate~$\pi/2$ phase ahead the depth~$h$.
That is because the strong turbulent mixing at the troughs makes the concentration peak quickly.

There is only one significant peak in one period in Figure~\ref{sedPRhu}--\ref{sedPRx}.
\cite{Zedler:2006fk}, who reported numerical results of large eddy simulation of the flow and suspended sediment over sinusoidal ripples, found three peaks on the time series of concentration at the crest and trough.
Their ripples  have a height to wavelength ratio of~$0.1$, which is five times steeper than our ripples.
 \cite{Zedler:2006fk} commented that these peaks are mainly due to the vortex, shear stress and advection near the ripple.
In our simulation, the only significant peak in a period is due to the turbulent mixing.

\begin{figure}
\centering
\begin{tabular}{c@{}c}
\rotatebox{90}{\hspace{15ex}depth~$h$} &
\includegraphics{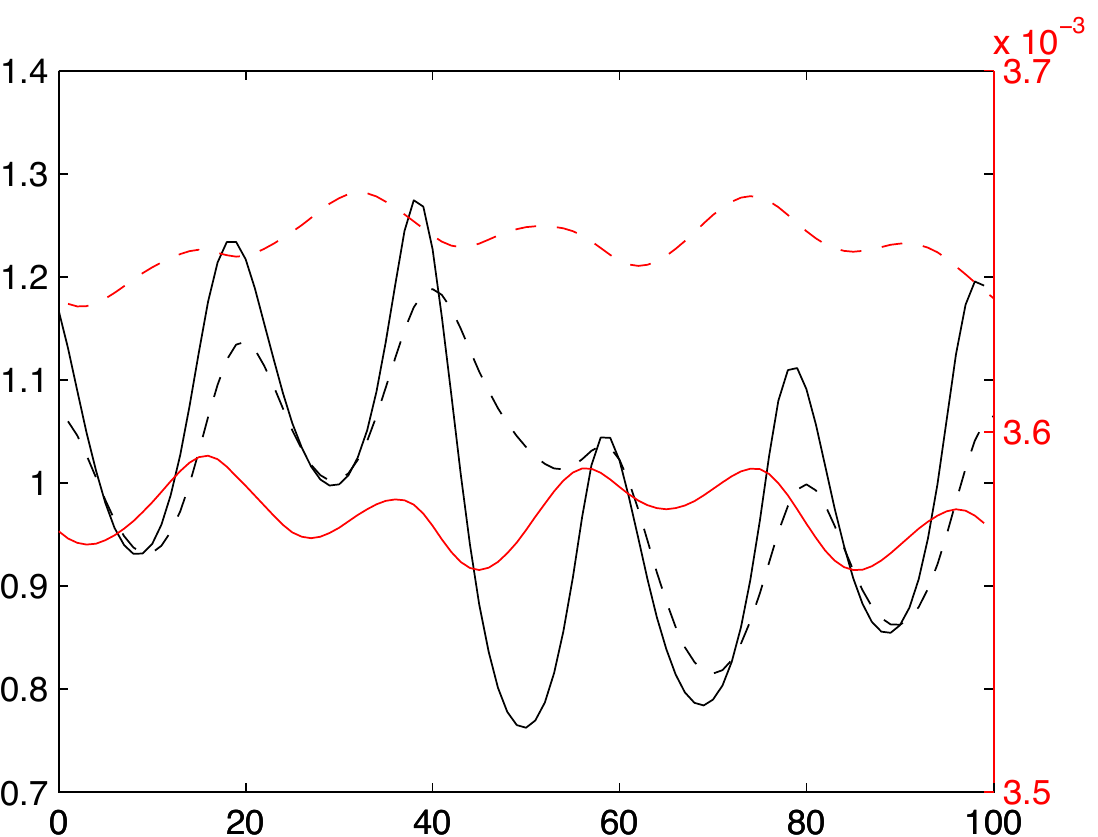}
\rotatebox{90}{\hspace{15ex}concentration~$\cc$}\\
& $x$
\end{tabular}
\caption{Plots of the depth~$h$~(black curves) and depth-averaged concentration~$\cc$~(red curves) of the fluid flowing over the rippled bed with ripple height~$0.4$~(dash curves) and~$0.6$~(line curves) in the $x$~direction at time $t=180$.
}
\label{sedPRsteep}
\end{figure}

Figure~\ref{sedPRsteep} compares at the time $t=180$ the depth~$h$~(black) and depth-averaged concentration~$\cc$~(red) of the flow over ripples with different heights.
The dash curves represent the depth~$h$ and depth-averaged concentration~$\cc$ for the ripple height~$0.4$, while the line curves are for the ripple height~$0.6$.
The fluid depth~$h$ is usually bigger at the crest of the steeper ripple, but the depth-averaged concentration~$\cc$ becomes smaller for the steeper ripples.
However, the laboratory experiment by~\cite{Osborne:1996fk}, whose ripples have an approximate height to wavelength ratio of~$0.2$, verifies that steep asymmetric ripples under shoaling waves produce greater concentrations higher in the water column than low steepness ripples.
Such difference is possibly because the ripples with small height in our simulation do not produce strong vortices that enhance the pick up of sediment into suspension.
The phenomenon of smaller depth-averaged concentration for steeper ripples is because the increased fluid depth for steeper ripples produces small depth-averaged concentration according to the erosion and deposition $w_f\cc/h$ in the governing equation~\eqref{sed:cH}.

\section{Conclusion}
\label{sec:con}

This work derives a suspended sediment model~\eqref{smag:h}--\eqref{sed:cH} to simulate the interactions between suspended sediment and turbulent flows.
The concentration equation~\eqref{sed:cH} consists of the effects of sediment erosion, advection, and dispersion.
Section~\ref{smag:Tbc} embedded the physical boundary conditions on the free surface and on the mean bed in a family of problems to access a slow manifold in the system.
The parameter~$\gamma=\gamma_c=1$ recovers the original physical problem.
Based on the small variations $\gamma=\gamma_c=\tan\theta=w_f=\partial_x=\partial_y=0$, a four parameter family of equilibria exists to support the existence of a slow manifold in the system, a slow manifold describing the large lateral structures.
Computer algebra detailed in Appendix~\ref{caSmag} leads to the evolution equations in the field of depth $h(x,y,t)$, depth-averaged lateral velocities $\uu(x,y,t)$ and $\vv(x,y,t)$, and depth-averaged concentration~$\cc(x,y,t)$.
It is reassuring that the dominant terms in our model \eqref{sed:cH} agree with  established modelling~\cite[e.g.]{Wu:2004uq,Duan2004,Duan:2006fk}.
Then our model includes more subtle effects, that could be important for suspended sediment in complex flow regimes.
The trends of the suspended sediment concentration corresponds to the published experimental measurements~\cite[e.g.]{Cellino:1999fk,Yoon:2005fk}.

Section~\ref{sec:Num} implemented numerical simulations of the suspended sediment under large waves by the suspended sediment model~\eqref{smag:h}--\eqref{sed:cH}.
The time series of the depth-averaged suspended sediment concentration rises fast and falls slowly.
Supercritical flow arises when the fluid flowing over the rippled bed.
The plots of the depth-averaged concentration~$\cc$ in space show that high concentration arises at the troughs and low concentration arises at the crests.

\bibliographystyle{agsm}
\bibliography{bibexport}

\appendix
\section{Ancillary computer algebra program}
\label{caSmag}

This appendix lists the computer algebra to construct the slow manifold model of the suspended sediment in turbulent flow.

Denote the fluid depth~$h(x,y,t)$ by~\verb|h|, depth-averaged lateral velocities $\uu(x,y,t)$ and~$\vv(x,y,t)$ by \verb|uu|~and~\verb|vv|, depth-averaged suspended sediment concentration $\cc(x,y,t)$ by~\verb|cc|,
and their time derivatives by $h_t=\verb|gh|$, $\uu_t=\verb|gu|$, $\vv_t=\verb|gv|$ and $\cc_t=\verb|gc|$.  
Denote the mean bed~$b(x,y)$ by~\verb|b|.
The coefficients of lateral and normal gravitational forcing are represented by~\verb|grx|, \verb|gry|~and~\verb|grz:=1|.  
Use~\verb|qq| represent the mean flow speed~$\bq=\sqrt{\uu^2+\vv^2}$ and~$\verb|rqq|$ for the reciprocal of this mean speed.

Use the operator \verb|h(m,n)| to denote the various lateral derivatives of the fluid depth,~$\partial_x^m\partial_y^nh$.
Similarly use the operators~\verb|uu(m,n)|,~\verb|vv(m,n)|, \verb|cc(m,n)| to denote lateral derivatives of the depth-averaged lateral velocities~$\uu(x,y,t)$ and~$\vv(x,y,t)$, and the depth-averaged concentration~$\cc(x,y,t)$.
These operators depend upon time and lateral space.
Then the lateral derivative $\partial_x\tt{h(m,n)}=\partial_x^{m+1}\partial_y^n\tt h$, and the time derivative~$\partial_t\tt{h(m,n)}=\partial_x^m\partial_y^n\tt{gh}$, for example.
Define readable abbreviations for~$h(x,y,t)$ and its first spatial derivatives.
We use~\verb|d| to count the number of lateral derivatives so we can easily truncate the asymptotic expansion.

Define the operators for the mean flow speed~$\bq$ and its reciprocal.
The last simplification rule for~{\tt rqq} breaks the symbolic symmetric between the two lateral directions.
However the benefit of canonical representation outweighs the cost of loss of symbolic symmetry.

The key to the correctness of this program is that the residuals of the governing equations are computed correctly, and that the algorithm only terminates when these residuals are zero to the specified error.

\small
\verbatimlisting{caSedFluid.txt}

\end{document}